\documentclass[10pt,fleqn]{article}
\usepackage[fleqn]{amsmath}
\usepackage{amssymb, latexsym}
\usepackage[utf8]{inputenc}
\usepackage[T2A]{fontenc}
\usepackage[english]{babel}
\usepackage{amsthm}
\usepackage{xcolor}
\usepackage{array}
\usepackage{proof}

\theoremstyle{nonumberplain}
\newtheorem{lem}{Lemma}
\newtheorem{theo}{Theorem}

\newtheorem{definition}{Definition}
\newtheorem{remark}{Remark}

\theoremstyle{plain}

\newtheorem{theorem}{Theorem}
\newtheorem{proposition}[theorem]{Proposition}
\newtheorem{lemma}[theorem]{Lemma}
\newtheorem{corollary}[theorem]{Corollary}
\newtheorem{fact}[theorem]{Fact}

\newcommand{\blem}{\begin{lem}}
\newcommand{\elem}{\end{lem}}
\newcommand{\btheo}{\begin{theo}}
\newcommand{\etheo}{\end{theo}}
\newcommand{\bl}{\begin{lemma}}
\newcommand{\el}{\end{lemma}}
\newcommand{\bt}{\begin{theorem}}
\newcommand{\et}{\end{theorem}}
\newcommand{\bcor}{\begin{corollary}}
\newcommand{\ecor}{\end{corollary}}
\newcommand{\bp}{\begin{proof}}
\newcommand{\ep}{\end{proof}}
\newcommand{\bpr}{\begin{proposition}}
\newcommand{\epr}{\end{proposition}}
\newcommand{\brem}{\begin{remark}}
\newcommand{\erem}{\end{remark}}
\newcommand{\bd}{\begin{definition}}
\newcommand{\ed}{\end{definition}}
\newcommand{\bex}{\begin{example} }
\newcommand{\eex}{\end{example}}
\newcommand{\bfc}{\begin{fact}}
\newcommand{\efc}{\end{fact}}

\newcommand{\beq}{\begin{equation*}}
\newcommand{\eeq}{\end{equation*}}
\newcommand{\bgas}{\begin{gather*}}
\newcommand{\egas}{\end{gather*}}
\newcommand{\ben}{\begin{enumerate}}
\newcommand{\een}{\end{enumerate}}

\newcommand{\model}{\vDash}
\newcommand{\Iff}{\Leftrightarrow}

\newcommand{\psaid}{p\:\mathtt{said}~}

\renewcommand{\and}{\wedge}
\newcommand{\se}{\vdash}
\renewcommand{\phi}{\varphi}
\newcommand{\imp}{\rightarrow}

\newcommand{\Imp}{\Rightarrow}
\newcommand{\Eqv}{\Leftrightarrow}

\begin{document}
\begin{center}\Large{\textbf{Proof Systems and Models for the\\ First-Order Primal Logic}}\\
\end{center}
\begin{center} Alexandra Podgaits
\end{center}
\begin{center} Moscow State University, Moscow, Russian Federation\\
podgaitsal@gmail.com
\end{center}
\begin{center}\begin{tabular}{p{10cm}}\textbf{Abstract.} We study the first-order primal infon logic. It is the core of the policy language DKAL. We provide Gentzen-style calculi for two versions of this logic that are not equivalent. For both versions we investigate the semantics: one of them is a generalization of the so-called quasi-boolean semantics, the other one is a Krypke-style semantics. We prove the completeness results and the disjunction property for both logics.
\end{tabular}
\end{center}

\section{Introduction}
In this paper we consider the first-order primal infon logic. Primal infon logic was introduced by Y. Gurevich and I. Neeman in \cite{simple} in connection with the policy language DKAL (Distributed Knowledge Authorization Language) [1-3, 9].

In the primal infon logic every statement is considered as a piece of information (infon) that is used in communication between agents. A conjunction of infons $\phi\and\psi$ is considered as the least information containing both $\phi$ and $\psi$, an implication $\phi\imp\psi$ is understood as the least information from which an agent can obtain $\psi$ as soon as it knows $\phi$.

For every principal $p$ there is an infon $\psaid\phi$. The intuitive meaning of $\psaid\phi$ is that $\phi$ can be derived from the information sent by $p$. In the infon logic, $\psaid\phi$ is considered as a modality satisfying the axioms of the modal logic K. These modalities are sometimes called quotation modalities.

So far, the propositional primal infon logic was mostly studied. In~\cite{simple} the infon logic was introduced as the $\{\and,\imp\}$-fragment of the intuitionistic logic extended by quotation modalities. For the reasons of efficiency the infon logic was restricted to the so-called \textit{primal} logic, a version of the intuitionistic logic with a weak form of implication. Y. Gurevich and I. Neeman in \cite{datalog} proved that the derivability problem for the propositional primal logic is linear-time decidable, and this result was extended by Y. Gurevich and C. Сotrini in \cite{cotrini} to the primal logic with a weak form of disjunction and unrestricted quotation modalities.

In~\cite{simple} a Kripke-style semantics for the primal infon logic was introduced. L. Beklemishev and Y. Gurevich \cite{prop} added disjunction to the language and introduced a simpler semantics for the primal infon logic - the so-called quasi-boolean semantics, which is a modification of the boolean semantics for the classical logic.

The logic implemented in DKAL is in fact a weak form of the first-order primal logic. 


In this paper we investigate two versions of the first-order primal logic. We would like to thank the referee of the previous version of this paper for finding a mistake in it. Understanding the cause of the mistake led to discovering that the two versions of the first-order primal logic differ from each other to a greater extent than we thought. As it turns out, the two types of semantics corresponding for the propositional primal logic behave differently in the first order case. A natural generalization of the intuitionistic Kripke-style semantics for the propositional primal logic corresponds to the extension of a single-conclusion Gentzen-style calculus by the standard quantifier rules. On the other hand, a natural generalization of the quasi-boolean semantics leads to a stronger logic. This logic can be alternatively axiomatized either by a similar extension of the multi-conclusion Gentzen-style calculus, or by adjoining to the single-conclusion calculus the well-known \emph{constant domain principle}: \begin{gather}
\forall x(A\vee B(x))\Imp A\vee\forall x B(x).
\end{gather}

For both versions of predicate primal logic we prove soundness and completeness theorems by adapting the methods from \cite{takeuti}.

The intuinionistic logic extended by (1) is sound and complete w.r.t. intuitionistic models with constant domains (\textit{Grzegorczyk's models}) (\cite{goernemann}). Unlike for the primal version, there is no known axiomatization by a simple cut-free Gentzen-style calculus for the intuitionistic logic of constant domains (\cite{fitting, mints, le}). The more complex axiomatizations of this logic can be found in \cite{nested, shimura}.

\section{First-Order Primal Logic}
Basic symbols of our language are: predicate symbols, constants, logical connectives $\top$, $\and$, $\vee$, $\imp$, $\forall$ and $\exists$  and disjoints sets of free and bound variables ($a,b,c,\dots$ and $x,y,z,\dots$, respectively). Formulas of the language are defined in the standard way.

Notice that we consider the language without function symbols.

\subsection{Natural Deduction Calculus }
In this section we consider natural deduction calculus in sequential format from \cite{troelstra}.

We define the relation $\Gamma\se\phi$, \textit{formula $\phi$ is derivable from the set of hypotheses $\Gamma$}, as the minimal relation containing the following axioms and closed with respect to the following inference rules:

Axioms:
\begin{gather*}
\se\top \quad\quad \phi\se\phi
\end{gather*}
Rules: 
\begin{gather*}
\begin{array}{ll}
\text{($\and$E)\quad} \dfrac{\Gamma\se\phi\and\psi}{\Gamma\se\phi} \quad \dfrac{\Gamma\se\phi\and\psi}{\Gamma\se\psi} \qquad\qquad\quad & \text{($\and$I)\quad}  \dfrac{\Gamma\se\phi\quad\Gamma\se\psi}{\Gamma\se\phi\and\psi}  \\
\text{($\vee$I)\quad} \dfrac{\Gamma\se\phi}{\Gamma\se\phi\vee\psi} \quad \dfrac{\Gamma\se\psi}{\Gamma\se\phi\vee\psi} & \text{($\vee$E)\quad}\dfrac{\Gamma,\theta\se\phi\quad\Gamma,\psi\se\phi\quad\Gamma\se\theta\vee\psi}{\Gamma\se\phi}\\
\text{($\imp$E)\quad}  \dfrac{\Gamma\se\phi\quad\Gamma\se\phi\imp\psi}{\Gamma\se\psi} & \text{($\imp$IW)\quad}  \dfrac{\Gamma\se\psi}{\Gamma\se\phi\imp\psi} \\
\text{($\forall$I)\quad} \dfrac{\Gamma\se \phi(a)}{\Gamma\se\forall x \phi(  x)} & \text{($\forall$E)\quad} \dfrac{\Gamma\se\forall  x\phi(  x)}{\Gamma \se \phi[  u/  x]} \\
\text{where} \quad   a\not\in FreeVar(\Gamma),&   u\; \text{is a constant or a variable} \\
\text{($\exists$I)\quad}  \dfrac{\Gamma\se\phi[  u/ x]}{\Gamma\se\exists x\phi( x)} & \text{($\exists$E)}\quad \dfrac{\Gamma\se\exists  x\phi( x) \quad\Delta, \phi( a)\se \psi}{\Gamma,\Delta\se \psi} \\
 u\; \text{is a constant or a variable,}  &  a\not\in FreeVar(\Gamma,\Delta) \\
\text{(Weakening)	\quad}   \dfrac{\Gamma\se\psi}{\Gamma,\phi\se\psi} & \text{(Trans)\quad}  \dfrac{\Gamma\se\phi\quad\Gamma,\phi\se\psi}{\Gamma\se\psi} \\
\end{array}
\end{gather*}

This system is denoted QP; its propositional fragment is denoted P.

Notice that the deduction theorem fails for primal logic; for example, the formula $\forall x\phi(x)\imp\forall x\phi(x)$ is not derivable in the primal logic.

We also consider primal logic with \textit{weak disjunction}, which is obtained from primal logic by removing the rule ($\vee$E). The corresponding systems are denoted QPW and PW.

\subsection{ Gentzen-Style Calculi}
Sequents are objects of the form $\Gamma\Imp\Delta$, where $\Gamma,\Delta$ are finite sets of formulas.
We consider two types of Gentzen-style calculus: multiple-conclusion format and version where $\Delta$ consists of exactly one formula. For reasons of brevity we formulate the multiple-conclusion variant; to obtain the second one all rules must be modified in such a way that the succedent consists of a single formula.

The axioms and rules of the first-order primal logic are standard except for the rule ($\imp$Rp).

Axioms:
$$
\phi\Imp\phi \quad \quad \Imp\top
$$
Rules:
$$
\begin{array}{clcl}
\text{($\and$L)\quad} & \dfrac{\Gamma,\phi,\psi\Imp\Delta}{\Gamma,\phi\and\psi\Imp\Delta}\qquad\qquad\quad&
\text{($\and$R)\quad} & \dfrac{\Gamma\Imp\Delta,\phi\quad\Gamma\Imp\Delta,\psi}{\Gamma\Imp\Delta,\phi\and\psi}  \\
\text{($\vee$L)\quad}&\dfrac{\Gamma,\phi\Imp\Delta \quad \Gamma,\psi\Imp\Delta}{\Gamma,\phi\vee\psi\Imp\Delta}\quad&
\text{($\vee$R)\quad}&\dfrac{\Gamma\Imp\Delta,\phi,\psi}{\Gamma\Imp\Delta,\phi\vee\psi} \\
\text{($\imp$L)\quad} & \dfrac{\Gamma,\psi\Imp\Delta\quad\Gamma\Imp\Delta,\phi}{\Gamma,\phi\imp\psi\Imp\Delta}\quad &
\text{($\imp$Rp)\quad} & \dfrac{\Gamma\Imp\psi,\Delta}{\Gamma\Imp\phi\imp\psi,\Delta}\\
\text{($\forall$R)\quad} &\dfrac{\Gamma\Imp \Delta,\phi( a)}{\Gamma\Imp\Delta,\forall x \phi( x)} & \text{($\forall$L)\quad}& \dfrac{\Gamma,\phi[ u/ x]\Imp\Delta}{\Gamma, \forall x\phi( x) \Imp\Delta} \\
 & \text{where} \quad  a\not\in FreeVar(\Gamma,\Delta), & &\;  u \text{ is a } \\[-10pt]
 & & & \text{constant or a variable}\\
\text{($\exists$R)\quad} & \dfrac{\Gamma\Imp\Delta,\phi[ u/ x]}{\Gamma\Imp\Delta,\exists x\phi( x)} & \text{($\exists$L)\quad}& \dfrac{\Gamma,\phi( a)\Imp\Delta}{\Gamma,\exists x\phi( x)\Imp\Delta}\\
& \text{with the same constraints.} & &  \\
\text{(Weakening)\quad}&\dfrac{\Gamma\Imp\Delta}{\Gamma,\Gamma_1\Imp\Delta,\Delta_1}\quad &
\text{(Cut)\quad} & \dfrac{\Gamma\Imp\Delta,\phi\quad\phi,\Gamma_1\Imp\Delta_1}{\Gamma,\Gamma_1\Imp\Delta,\Delta_1} \\
\end{array}
$$
We denote this calculus QGPM; the corresponding propositional calculus is denoted GPM.

To obtain a system with weak disjunction one must remove the rule ($\vee$L). The corresponding systems with weak disjunction are denoted QGPMW and GPMW. 

If we let the succedent consist of exactly one formula and replace the rule ($\imp$L) with the following rule:
$$\text{($\imp$L)\quad} \dfrac{\Gamma,\psi\Imp\theta\quad\Gamma\Imp\phi}{\Gamma,\phi\imp\psi\Imp\theta},$$
and the rule ($\vee$R) with the following rule:
$$ \text{($\vee$R)\quad} \dfrac{\Gamma\Imp\phi_i}{\Gamma\Imp\phi_1\vee\phi_2}, i=1,2,$$
we obtain the calculus that will be denoted QGP (GP for the propositional variant, GPW and QGPW for systems with weak disjunction).

Unlike for the intuitionistic version of QGPM \cite{fitting}, for both QGPM and QGP the cut-rule can be eliminated.
\bt \label{cutm}
A sequent $\Gamma\Imp\Delta$ is derivable in \emph{QGPM} iff it is derivable in \emph{QGPM} without the use of \emph{(Cut)}.
\et
There are both syntactical and semantical proofs for Theorem \ref{cutm}. A syntactical proof can be obtained by the standard techniques (for example, from \cite{takeuti}). The only new reductions concern the weak implication rules and they were considered in \cite{prop}. A semantical proof will be given in section 2.

\bt \label{cut}
A sequent $\Gamma\Imp\Delta$ is derivable in \emph{QGP} iff it is derivable in \emph{QGP} without the use of \emph{(Cut)}.
\et
A syntactical proof of this theorem can be obtained by the standard techniques (for example, from \cite{takeuti}).

For the propositional case two formulated Gentzen-style calculi are equivalent.
\bpr
A sequent $\Gamma\Imp\bigvee\Delta$ is derivable in \emph{GP} iff $\Gamma\Imp\Delta$ is derivable in \emph{GPM}.
\epr
For the first-order case only one implication holds:
\bpr
If a sequent $\Gamma\Imp\bigvee\Delta$ is derivable in \emph{QGP} then $\Gamma\Imp\Delta$ is derivable in \emph{QGPM}.
\epr
The converse is not true: we provide a counterexample.
Recall that \emph{constant domain principle} is a schema
$$\text{(CD)}\quad\forall x(A\vee B(x))\Imp A\vee\forall x B(x).$$
\bl\label{constdomain}
Constant domain principle is derivable in \emph{QGPM} but is not derivable in \emph{QGP}.
\el
\textbf{Proof.}
First we provide a derivation of constant domain principle in QGPM.
\beq
\dfrac{ \dfrac{  \dfrac{\dfrac{    \dfrac{A\Imp A}{A\Imp A,B(a)}\qquad \dfrac{B(a)\Imp B(a)}{B(a)\Imp A,B(a)  }}{ A\vee B(a)\Imp A,B(a)}}{\forall x(A\vee B(x))\Imp A, B(a)}  }{  \forall x(A\vee B(x))\Imp A,\forall x B(x) }  }{\forall x(A\vee B(x))\Imp A\vee\forall x B(x)}
\eeq

Suppose that CD is derivable in QGP without (Cut).

Consider the lowermost rule application in the derivation. It could be ($\forall$L) or ($\vee$R).
\begin{enumerate}
\item The lowermost rule application was ($\forall$L):
\begin{gather*}
\dfrac{A\vee B(u)\Imp A\vee\forall x B(x)}{\forall x(A\vee B(x))\Imp A\vee\forall x B(x)}
\end{gather*}

Two cases are possible:

\begin{enumerate}
\item Before ($\forall$L) ($\vee$L) was applied:
\begin{gather*}
\dfrac{A\Imp A\vee\forall x B(x) \quad B(u)\Imp A\vee\forall x B(x) }{A\vee B(u)\Imp A\vee\forall x B(x)}.
\end{gather*}
The right sequent is not derivable, for $B(u)$ implies neither $A$ nor $\forall x B(x)$.
So this case is impossible.

\item Before ($\forall$L) ($\vee$R) was applied:
\begin{gather*}
\dfrac{A\vee B(u)\Imp A}{A\vee B(u)\Imp A\vee\forall x B(x)}\quad\mbox{or} \quad\dfrac{A\vee B(u)\Imp \forall x B(x)}{A\vee B(u)\Imp A\vee\forall x B(x)}
\end{gather*}
But $A\vee B(u)$ implies neither $A$ nor $\forall x B(x)$. Contradiction.
\end{enumerate}

\item The lowermost rule application was ($\vee$R):
\ben
\item
\beq
\dfrac{\forall x(A\vee B(x))\Imp \forall x B(x)}{\forall x(A\vee B(x))\Imp A\vee\forall x B(x)}
\eeq
Before ($\vee$R) could be ($\forall$L) or ($\forall$R):
\beq
\dfrac{A\vee B(u)\Imp \forall x B(x)}{\forall x(A\vee B(x))\Imp\forall x B(x)}\quad\mbox{or} \quad \dfrac{\forall x(A\vee B(x))\Imp B(a)}{\forall x(A\vee B(x))\Imp\forall x B(x)}
\eeq
It can be easily shown that neither can be the case.
\item
\beq
\dfrac{\forall x(A\vee B(x))\Imp A}{\forall x(A\vee B(x))\Imp A\vee\forall x B(x)}
\eeq
In this case before ($\vee$R) ($\forall$L) was applied:
\beq
\dfrac{A\vee B(u)\Imp A  }{\forall x(A\vee B(x))\Imp A}
\eeq
This cannot be the case.
\een
Contradiction.\qed
\end{enumerate}
Actually, if we add constant domain principle to QGP we obtain exactly QGPM:
\bt
$\text{\emph{QGPM}}\se\Gamma\Imp\Delta\;$ iff $\;\text{\emph{QGP+(CD)}}\se\Gamma\Imp\bigvee\Delta$.
\et
\textbf{Proof.}
The implication form right to left is obvious.

We prove the converse by induction on the derivation of $\Gamma\Imp\Delta$.

Consider the lowermost rule application in the given derivation. The rules ($\and$L), ($\vee$L), ($\forall$L) and ($\exists$L) are the same in QGPM and QGP.

Consider the case when the lowermost rule application was ($\forall$R) (actually, this is the only interesting case):

$$ \dfrac{\Gamma\Imp\Delta,\phi(a)}{\Gamma\Imp\Delta,\forall x\phi(x)}, $$

where $ a\not\in FreeVar(\Gamma,\Delta)$.

Let $\delta = \bigvee\Delta$.

By induction hypothesis $QGP\se \Gamma\Imp\delta\vee\phi(a) $. We have to prove that QGP$\;\se\Gamma\Imp\delta\vee\forall x\phi(x)$.
$$
\infer
{\Gamma\Imp\delta\vee\forall x\phi(x) }
{\infer
{\Gamma\Imp\forall x (\delta\vee\phi(x))} 
{\Gamma\Imp\delta\vee\phi(a)} 
\quad \forall x (\delta\vee\phi(x))\Imp \delta\vee\forall x \phi(x)\;\mbox{(CD)}
}
$$

Consider the case when the lowermost rule application was ($\imp$L):
$$\dfrac{\Gamma,\psi\Imp\Delta\quad\Gamma\Imp\Delta,\phi}{\Gamma,\phi\imp\psi\Imp\Delta} $$

Let $\delta = \bigvee\Delta$.

By induction hypothesis $QGP\se \Gamma,\psi\Imp\delta$ and $QGP\se\Gamma\Imp\delta\vee\phi$. We have to prove that $QGP\se\Gamma,\phi\imp\psi\Imp\delta$.
$$
\infer
{\Gamma,\phi\imp\psi\Imp\delta}
{\infer{\Gamma,(\delta\vee\phi)\imp\psi\Imp\delta}{\Gamma,\psi\Imp\delta\quad\Gamma\Imp\delta\vee\phi}   \quad \infer
{\Gamma,(\delta\vee\phi)\imp\psi\Imp\phi\imp\psi}
{ \infer{\Gamma,\psi\Imp\phi\imp\psi}{\Gamma,\psi\Imp\psi} \quad \Gamma\Imp\delta\vee\phi   }}
$$
  
Other cases are similar.
\qed
  
It turns out that QGP is equivalent to the natural deduction calculus QP and QGPM is not (since QGPM and QGP are not equivalent).
\bt
$\Gamma\se\phi$ is derivable in \emph{QP} iff $\;\Gamma\Imp\phi$ is derivable in \emph{QGP}.
\et
\textbf{Proof.}
We prove both implications by induction on the derivation.

Suggest QGP$\;\se\Gamma\Imp\phi$. We consider the case when the lowermost application of a derivation rule was ($\exists$L); other cases are similar or trivial.

$$\dfrac{\Gamma,\phi( a)\Imp\psi}{\Gamma,\exists x\phi( x)\Imp\psi}$$

By induction hypothesis $\Gamma,\phi( a)\se\psi$ is derivable in QP.
So we can derive
$$\text{($\exists$E)}\quad \dfrac{\exists  x\phi( x)\se \exists  x\phi( x)\quad \Gamma,\phi( a)\se\psi}{\Gamma,\exists x\phi( x)\se\psi}$$

Now suggest $\Gamma\se\phi$ in QP. Consider the case when the lowermost rule application was ($\exists$E); other cases are similar or trivial.
$$\dfrac{\Gamma\se\exists  x\phi( x) \quad\Delta, \phi( a)\se \psi}{\Gamma,\Delta\se \psi}$$
By induction hypothesis $\Gamma\Imp\exists  x\phi( x)$ and $\Delta, \phi( a)\Imp \psi$ are derivable in QGP.
By applying ($\exists$L) and (Cut), we obtain $\Gamma,\Delta\Imp \psi$.
\qed

\section{Semantics}\label{semantics}
Y.Gurevich, I.Neeman \cite{first} introduced a Kripke-style semantics for the primal logic without disjunction. It is based on the notion of cone in a Kripke model and in fact is a modification of the standard Kripke semantics for the intuitionistic logic. A simpler semantics for the propositional primal logic with or without disjunction, the so-called \textit{quasi-boolean models}, was introduced in \cite{prop}. It can be considered as a non-deterministic modification of boolean valuations for the classical logic.

In this section we consider a natural generalization of quasi-boolean models to the first-order logic. We show that QGPM is sound and complete with respect to this semantics. However, QGP turns out to be sound but incomplete with respect to the quasi-boolean semantics. We introduce an extension of Kripke models from \cite{first} to the first-order logic and show that QGP is sound and complete with respect to it.  

\subsection{First-Order Quasi-boolean Semantics}
%
%


In this section we provide a natural generalization of quasi-boolean models to the first-order logic. For the definitions of quasi-boolean semantics in the propositional case see \cite{prop}.

A \textit{first-order quasi-boolean model} is a pair $\mathcal{M}=(M,v)$, where $v$ is a \textit{quasi-boolean valuation} that assigns  $0$ or $1$ to all closed atomic formulas and to all formulas of the form $\phi\imp\psi$, $M\ne\varnothing$, for every constant $c$ its interpretation $\tilde c\in M$ is specified. 

The forcing relation $M\model_v$ (we will denote it $\mathcal M\model$) naturally extends the valuation to all formulas:
\begin{enumerate}
\item $\mathcal M\model \phi \Iff  v(\phi)=1$, where $\phi$ is a closed atomic formula or an implication;
\item $\mathcal M\model\phi\and\psi\Eqv\quad\mathcal M\model\phi\text{  and }\mathcal M\model\psi$;
\item $\mathcal M\model\phi\vee\psi\Eqv\quad\mathcal M\model\phi\text{  or }\mathcal M\model\psi$;
\item $\mathcal M\model \forall x \phi( x) \Iff \forall u\in M \; \mathcal M\model \phi[ u/ x]$;
\item $\mathcal M\model \exists x \phi( x) \Iff \exists u\in M \; \mathcal M\model \phi[ u/ x]$.
\end{enumerate}

We say that a valuation $v$ is \textit{quasi-boolean}, if the following conditions hold for every implication $\phi\imp\psi$:
\ben
\item If $\mathcal M\model \psi$, then $\mathcal M\model \phi\imp\psi$;
\item If $\mathcal M\model \phi\imp\psi$, then either $\mathcal M\not\model \phi$ or $\mathcal M\model \psi$.
\een

We add a new variable $\bar{c}$ for every $c\in M$ to our language. We say that  $\Gamma(\vec a)\Imp\Delta(\vec a)$, where $\vec a$ is a vector of free variables, is valid in $\mathcal M$ iff for all $\vec u\in M $
$$\mathcal M\model\bigwedge\Gamma[\vec{\bar u}/ \vec a] \;\;\mbox{implies} \;\;\mathcal M\model\bigvee\Delta[\vec{\bar u}/ \vec a].$$ (We consider the valuation of the empty-set be $0$.)

QGPM is sound and complete with respect to the first-order quasi-boolean semantics. 
\bt\label{sound}
If $\Gamma\Imp\Delta$ is derivable in \emph{QGPM} then $\Gamma\Imp\Delta$ is valid in every first-order quasi-boolean model.
\et
Now we prove the completeness adapting the method from~\cite{takeuti}.
\bt\label{foquasi}
If $\Gamma\Imp\Delta$ is not derivable in \emph{QGPM} without the use of \emph{(Cut)} then there is a first-order quasi-boolean model such that $\Gamma\Imp\Delta$ is not valid in it. 
\et
\textbf{Proof.} See appendix.
$$$$
Now we can prove Theorem~\ref{cutm}.\\
\textbf{Corollary.}\emph{
A sequent $\Gamma\Imp\Delta$ is derivable in \emph{QGPM} iff it is derivable in \emph{QGPM} without the use of \emph{(Cut)}.}\\
\textbf{Proof.}
Consider $\Gamma\Imp\Delta$ that is not derivable in QGPM without the use of (Cut). By theorem~\ref{foquasi}, there is a first-order quasi-boolean model $\mathcal M$ such that $\Gamma\Imp\Delta$ is not valid in it. By theorem~\ref{sound}, $\Gamma\Imp\Delta$ is not derivable in QGPM (if it was, it would be valid in $\mathcal M$).
\qed

\textbf{Remark.} It is easy to prove that QGP is sound with respect to the first-order quasi-boolean models. However, the completeness does not hold: the constant domain principle $$\text{(CD)}\quad\forall x(A\vee B(x))\Imp A\vee\forall x B(x)$$ is valid in every first-order quasi-boolean model, but is not derivable in QGP.

\subsection{Semantics for QGP}\label{firstsem}
In this section we introduce Kripke models for the single-conclusion variant of the first-order primal logic. These are Kripke models for the first-order intuitionistic logic with a special condition on forcing for implications (and for disjunctions, is case of the logic with weak disjunction).\\

A \textit{Kripke model} for the first-order primal logic QGPM is a tuple \\$\mathcal W = \langle W, \leqslant, D, \nu\rangle$ such that
\ben
\item $(W,\leqslant)$ is a non-empty partially ordered set;
\item $D$ is a function assigning non-empty sets to the elements of $W$ such that $\forall u,v \in W (u\leqslant v \Imp D(u)\subseteq D(v))$;
\item $\nu$ is a map assigning $0$ and $1$ to tuples $\langle u, \phi(c_1,\dots, c_n)\rangle$, $u\in W, \tilde{c_i}\in D(u), \; \phi$ is a closed atomic formula or an implication; $\forall u,v \in W (u\leqslant v \Imp (\nu(\langle u, \phi(c_1,\dots, c_n)\rangle = 1 \Imp \nu(\langle v, \phi(c_1,\dots, c_n)\rangle = 1)$.
\een
The relation $\mathcal W, u \model \phi$ is defined as follows:
\begin{enumerate}
\item $\mathcal W, u \model \phi$ iff $\nu(u,\phi)=1$ for closed atomic formulas and implications;
\item $\mathcal W, u \model \top$;
\item $\mathcal W, u \model \phi\and\psi \Iff (\mathcal W, u \model \phi \text{  and  }\mathcal W, u \model \psi)$;
\item $\mathcal W, u \model \phi\vee\psi \Iff (\mathcal W, u \model \phi \text{  or  }\mathcal W, u \model \psi)$;
\item $\mathcal W, u\model \exists x \phi( x) \Iff (\exists c\in D(u) \; (\mathcal W, u\model \phi[ c/ x]))$;
\item $\mathcal W, u\model\forall x \phi( x) \Iff (\forall v\geqslant u \;( \forall c\in D(v) \;\mathcal W, v\model \phi[ c/ x])).$
\end{enumerate}
Also $\mathcal W, u \model \phi$ should satisfy the following conditions:
\ben
\item $\mathcal W, u \model \psi \Imp \mathcal W, u \model \phi\imp\psi $;
\item $\mathcal W, u \model \phi\imp\psi \Imp \forall v \geqslant u \;(\mathcal W, v \not\model \phi \text{  or  } \mathcal W, v \model \psi)$.
\een

We add a new variable for every constant $c \in \bigcup_{u\in W}D(u)$ and denote it $\bar c$.

We say that $\Gamma( \vec a)\Imp\phi( \vec a)$ is valid in $\mathcal W$ iff 
$$\forall u\in W\;\forall \vec c\in D(u) \;  (\mathcal W,u\model\bigwedge\Gamma( \vec{\bar c})\Imp\mathcal W,u\model\phi( \vec{\bar c})).$$

\textbf{Remark.}
To obtain a Kripke model for the first-order primal logic with weak disjunction (a \emph{WD-model}) replace 4. with
$$4'. \;\mathcal W, u \model \phi_i\Imp \mathcal W, u \model \phi_1\vee \phi_2, i=1,2.$$

QGP is sound and complete with respect to Kripke models. The proof of soundness is routine; to prove completeness it is convenient to consider natural deduction calculus.
\bt[Soundness] If a sequent $\Gamma\Imp\phi$ is derivable in \emph{QGP} then $\Gamma\Imp\phi$ is valid in every Kripke model $\mathcal W$.
\et

\bt[Completeness]\label{kripke} If a sequent $\Gamma\se\phi$ is valid in every Kripke model, then $\Gamma\se\phi$ is derivable in \emph{QP}.
\et
\textbf{Proof.}
The proof is similar to the standard proof of completeness for intuitionistic logic.
We prove that if $\Gamma\not\se\phi$ then there is a Kripke model $\mathcal W$ such that $\mathcal W\not\model\Gamma\se\phi$.

Let $C$ be a set of constants. A set of sentences $\Gamma$ is \textit{$C$-saturated}, iff
\ben
\item $\Gamma\se\psi$ implies $\psi\in\Gamma$;
\item $\Gamma\se\alpha\vee\beta$ implies $\alpha\in\Gamma\text{  or  }\beta\in\Gamma$;
\item If $\Gamma\se\exists x\psi(x)$ then $\text{ for some }  c\in C\; \psi(c)\in\Gamma$.
\een
\bl\label{sat} Suppose $\Gamma\not\se\phi$, $\Gamma,\phi$ in a language $L$. Let $C=\{c_0,c_1,c_2,\dots\}$ be a countable set of constants not in $L$, and let $L(C)$ be $L$ extended with $C$. Then there is a $C$-saturated $\Gamma'$ such that $\Gamma\subseteq\Gamma'$ and $\Gamma'\not\se\phi$.
\el
\textbf{Proof.}
Let $\psi_n$ be the numeration of all formulas in the extended language $L(C)$. We construct the sequence of sets of formulas $\Gamma_n$:
\ben
\item $\Gamma_0 = \Gamma$.
\item For $n>0$ several cases are possible:
\ben
\item $\psi_n\not\in\Gamma_n$ and $\Gamma_n,\psi_n\not\se\phi$. Then $\Gamma_{n+1} = \Gamma_n\cup\{\psi_n\}$.
\item $\psi_n\in\Gamma$ and $\psi_n$ has the form $ \theta_1\vee\theta_2$. For some $i\in\{1,2\}$ $\Gamma_n,\theta_i\not\se\phi$ (otherwise $\Gamma_n\se\phi$). Fix this $i$ and let $\Gamma_{n+1} = \Gamma_n\cup\{\theta_i\}$.
\item $\psi_n\in\Gamma$ and $\psi_n$ has the form $\exists x\theta(x)$. Let $c_i$ be the first constant from $C$ such that $c_i\not\in FreeVar(\Gamma_n)$. Let $\Gamma_{n+1} = \Gamma_n\cup\{\theta(c_i)\}$.
\item Otherwise $\Gamma_{n+1} = \Gamma_n$.
\een
\een \qed

Now we define canonical model for the primal logic.

Let $C_0,C_1,\dots$ be countable disjoint sets of constants: $C_i=\{c_{i0},c_{i1},\dots\}$, all $c_{ij}\not\in\mathcal L$. We write $C^n$ for $C_0\cup C_1\cup\cdots\cup C_n$.

Canonical model $\mathcal W$ is a tuple $\langle W, \leqslant, D, \nu\rangle$, where
\ben
\item $W$ is a set of $\Gamma'\supseteq\Gamma$ such that $\mathcal L(\Gamma')=\mathcal L\cup C^n$ for some $n$ and $\Gamma'$ is $C^n$-saturated;
\item $\Gamma'\leqslant\tilde\Gamma \Iff \Gamma'\subseteq\tilde\Gamma $;
\item if $\mathcal L(\Gamma')=\mathcal L\cup C^n$ and $\Gamma'$ is $C^n$-saturated then $D(\Gamma')=C^n$;
\item $\nu(\psi)=1 \Iff \psi\in\Gamma$ for closed atomic formulas and implications.
\een
\bl For every $\Gamma'\in W$ and any sentence $\psi$ of $\mathcal L(D(\Gamma'))$\\ $\Gamma'\model\psi\;\Iff\;\psi\in\Gamma'$.
\el
\textbf{Proof.}
We prove this by induction on $\psi$.

If $\psi$ is atomic or an implication, the condition holds by definition.

For $\psi=\alpha\and\beta$ it is obvious.

If $\psi=\alpha\vee\beta$, then $\Gamma'\model\alpha\vee\beta\Iff (\Gamma'\model\alpha \text{  or  } \Gamma'\model\beta)$ $\Iff
(\alpha\in\Gamma'\text{  or  }\beta\in\Gamma') \Iff ((\alpha\vee\beta)\in\Gamma')$ since $\Gamma'$ is D($\Gamma'$)-saturated.

For $\psi=\exists x\beta(x)$ the condition holds by saturation properties.

Let $\psi=\forall x\beta(x)$. Suppose $\Gamma'\model\psi$ but $\psi\not\in\Gamma'$ ($\Gamma'$ is $C^n$-saturated). For any $c\in C_{n+1}$,
$\Gamma'\not\se\beta(c)$. By Lemma \ref{sat} there is a $C^{n+1}$-saturated $\tilde\Gamma\supseteq\Gamma'$ such that $\tilde\Gamma\not\se\beta(c)$. Hence $\tilde\Gamma\not\model\beta(c)$. Contradiction with $\Gamma'\model\forall x\beta(x)$.

If $\psi\in\Gamma'$, it is trivial to show that $\Gamma'\model\psi$.
\qed
\bl
The forcing relation $\Gamma'\model\psi$ is quasi-boolean.
\el
\textbf{Proof.}
If $\Gamma'\model\beta$ then $\beta\in\Gamma'$.  Hence $\alpha\imp\beta\in\Gamma'\Iff\Gamma'\model\alpha\imp\beta$ by definition.

If $\Gamma'\model\alpha\imp\beta$ then $\forall\tilde\Gamma\geqslant\Gamma'\;\tilde\Gamma\model\alpha\imp\beta$ by monotonicity.
If now $\tilde\Gamma\model\alpha$, then both $\alpha\and\beta$ and $\alpha$ are in $\tilde\Gamma$, hence $\beta\in\tilde\Gamma$ and so $\tilde\Gamma\model\beta$.
\qed

Now we finish the proof of completeness. For our $\Gamma\not\se\phi$ there is a $C^n$-saturated $\Gamma'$ such that $\Gamma\subseteq\Gamma'$ and $\Gamma\not\se\phi$. So in our model $\Gamma'\not\model\phi$.

Note that we can make a rooted model from ours: we have to leave only the $\Gamma'$-cone, that is, all $\tilde\Gamma\geqslant\Gamma'$. This cone is still a Kripke model.
\qed

\textbf{Remark.}
To prove that QPW is complete with respect to WD-models one has to make a little change in the proof of the previous theorem:
remove condition 3 from the definition of a $C$-saturated set and change the definition of valuation $v$ in the canonical model: $v(\psi)=1\Iff\psi\in\Gamma'$ for closed atomic formulas, implications and disjunctions.

\subsection{Disjunction Property}
In this section we prove the disjunction property for both QGP and QGPM.

Recall that the set of Harrop formulas is defined by the following grammar:
$$H = \top | A | H\and H | B\imp H | \forall x H(x), $$
where A is an atomic formula, B is a formula.

\bt[Disjunction property] Let L be any of the logics \emph{QGP} or \emph{QGPM}. If $\Gamma$ is a set of Harrop formulas and \emph{L}$\;\se\Gamma\Imp\alpha\vee\beta$, then \emph{L}$\;\se\Gamma\Imp\alpha$ or \emph{L}$\;\se\Gamma\Imp\beta$. 
\et

\textbf{Proof.}
For any set of formulas $\Gamma$ we define a first-order quasi-boolean model $\mathcal M = (M,v)$ as follows:
\ben
\item $M=Const(\Gamma)\cup Const(\alpha)\cup Const(\beta) \cup \{a\}$, where $a$ is a new constant;
\item If $\phi$ is atomic then $v(\phi) = 1 \Leftrightarrow QGPM\se \Gamma\Imp \phi$;
\item $v(\phi\imp\psi) = 1 \Leftrightarrow ($\emph{L}$\;\se\Gamma\Imp(\phi\imp\psi)\quad\text{and}\quad(\mathcal M\not\model\phi \;\text{or}\;\mathcal M\model\psi))$.
\een
\bl\label{harrop} If $\mathcal M\model\phi$ then \emph{L}$\;\se\Gamma\Imp\phi$.
\el
\textbf{Proof.}
We prove this by induction on $\phi$. The only interesting case is $\phi = \forall x\psi(x)$.

Let $\mathcal M\model\forall x\psi(x)$. Then $\forall c\in M\;\mathcal M\model\psi(c)$, so $\mathcal M\model\psi(a)$. By induction hypothesis \emph{L}$\;\se\Gamma\Imp\psi(a)$. Since $a$ is a new constant, \emph{L}$\;\se\Gamma\Imp\forall x\psi(x)$.\qed

\bl Valuation $v$ is quasi-boolean.
\el
\textbf{Proof.}
If $\mathcal M\model\psi$, then \emph{L}$\;\se\Gamma\Imp\psi$. Then \emph{L}$\;\se\Gamma\Imp\phi\imp\psi$ and so the right-hand side of condition 3 holds. So by condition 3 $\;\mathcal M\model\phi\imp\psi$.\\
If $\mathcal M\model\phi\imp\psi$, then by condition 3 $\;\mathcal M\not\model\phi$ or $\mathcal M\model\psi$.\qed

\bl\label{harrop2} If $\phi$ is Harrop, then $\mathcal M\model\phi\Leftrightarrow$ \emph{L}$\;\se\Gamma\Imp\phi$. 
\el 
\textbf{Proof.}
The implication from left to right holds by lemma~\ref{harrop}. 

We prove the converse by induction on $\phi$. We consider the cases when $\phi$ is an implication or $\forall x\psi(x)$; other cases are trivial.

Let $\phi=\theta\imp\psi$, $\psi\in H$, \emph{L}$\;\se\Gamma\Imp\theta\imp\psi$. We have to show that $\mathcal M\not\model\theta$ or $\mathcal M\model\psi$. Suppose $\mathcal M\model\theta$. Then \emph{L}$\;\se\Gamma\Imp\theta$, so \emph{L}$\;\se\Gamma\Imp\psi$. Since $\psi$ is Harrop, by the induction hypothesis $\mathcal M\model\psi$.

Let $\phi=\forall x\psi(x)$, $\psi\in H$, \emph{L}$\;\se\Gamma\Imp\forall x\psi(x)$. For every $c\in M\;$ \emph{L}$\;\se\Gamma\Imp\psi(c) $. By the induction hypothesis $\forall c\in M\;\mathcal M\model\psi(c)$, so $\mathcal M\model\forall x\psi(x)$.
\qed
Now we finish the proof of the disjunction property. Recall that both QGP and QGPM are sound with respect to first-order quasi-boolean models.

By lemma~\ref{harrop2} if $\theta\in\Gamma$ then $\mathcal M\model\theta$. Since \emph{L}$\;\se\alpha\vee\beta$, by soundness $\mathcal M\model\alpha\vee\beta$. Hence $\mathcal M\model\alpha$ or $\mathcal M\model\beta$. By lemma~\ref{harrop} \emph{L}$\;\se\Gamma\Imp\alpha$ or \emph{L}$\;\se\Gamma\Imp\beta$. \qed

\section{Conclusion}
We considered two first-order variants of the primal logic and introduced semantics for both of them. We proved the cut-elimination and completeness for both logics. However, the weaker one, QGP, seems to be more appropriate for the primal logic. The stronger one, QGPM, corresponds to an unlikely situation when all the agents have the same domain and know about that. 

From the semantical point of view QGP is more like intuitionistic logic when QGPM is more like classical logic. However, the disjunction property holds for both of them.

\newpage
\section*{Appendix}

We give a sketch of a syntactical proof of Theorem~\ref{cut}.\\
\textbf{Theorem \ref{cut}.}
\emph{A sequent $\Gamma\Imp\Delta$ is derivable in \emph{QGP} iff it is derivable in \emph{QGP} without the use of \emph{(Cut)}.}

\textbf{Proof.}
We follow the standard procedure due to Gentzen.

Consider an application of a cut-rule $d$. We define its grade $g(d) = |\theta| + 1$, where $\theta$ is the cut-formula in $d$, and $|\theta|$ is the height of its parse tree. The rank $r(d)$ is the sum of heights of the right and left proof subtrees of the proof-tree ending with $d$. The proof goes by induction on the grade and a subsidiary induction on the rank of $d$. 

Consider the cut-rule with maximum grade. If there are several such rules, take the one with maximum rank and the leftmost of them.

It the cut-formula was not introduced on both sides of the proof tree immediately before the cut, it is easy to obtain a proof tree where this cut will have lower rank.

Consider the case when the cut-formula was introduced on both sides of the proof tree immediately before the cut. We consider the only special case when the cut-formula has the form $\alpha\imp\beta$; other cases are standard.

$$
\infer
{\Gamma,\Gamma_1\Imp\theta}
{\infer{\Gamma\Imp\alpha\imp\beta}{\Gamma\Imp\beta} \quad
 \infer{\Gamma_1,\alpha\imp\beta\Imp\theta}{\Gamma_1,\beta\Imp\theta\quad\Gamma\Imp\alpha}
}
$$

We reduce it to a proof
$$
\infer{\Gamma,\Gamma_1\Imp\theta}
{\Gamma\Imp\beta\quad\Gamma_1,\beta\Imp\theta}
$$

with a cut-rule with lower rank and lower grade.
\qed

$$$$
Now we prove theorem \ref{foquasi}.\\
\textbf{Theorem \ref{foquasi}.}
\emph{If $\Gamma\Imp\Delta$ is not derivable in \emph{QGPM} without the use of \emph{(Cut)} then there is a first-order quasi-boolean model such that $\Gamma\Imp\Delta$ is not valid in it. }

\textbf{Proof.} 
For the proof of this theorem we need the notion of the reduction tree. A sequent is written in every node of the tree.
We say that a free variable is available if it appears in some sequent on some step.
We define the process of building the reduction tree of $\Gamma_0\Imp\Delta_0$ as follows:
\ben
\item[Step] 0. We write $\Gamma_0\Imp\Delta_0$ in the root of the tree.
\item[Step] $k$, $k>0$. Two cases are possible:
\ben
\item In every leaf of the tree antecedent and succedent have a common formula. The process is complete.
\item Let $\Gamma\Imp\Delta$ be a sequent written in some leaf of the tree such that $\Gamma\cap\Delta=\varnothing$. 11 cases are possible:
\een
\een
\ben
\item [$k\equiv$]$0\;\;$(mod 11)

Let $\phi_1\and\psi_1,\dots,\phi_n\and\psi_n$ be all the conjunctions in $\Gamma$ such that reductions have not been applied to them on the previous steps. We add a node above $\Gamma\Imp\Delta$ and write there the sequent
$$\phi_1,\dots\phi_n,\Gamma\Imp\psi_1,\dots,\psi_n,\Delta$$
We say that a ($\and$L) reduction was applied to $\phi_1\and\psi_1,\dots,\phi_n\and\psi_n$.
\item [$k\equiv$]$1\;\;$(mod 11)

 Let $\phi_1\and\psi_1,\dots,\phi_n\and\psi_n$ be all the conjunctions in $\Delta$ such that reductions have not been applied to them on the previous steps.  We add $2^n$ nodes above $\Gamma\Imp\Delta$ and write there all the sequents
$$\Gamma\Imp\theta_1,\dots,\theta_n,$$
where $\theta_i$ is either $\phi_i$ or $\psi_i$.

We say that a ($\and$R) reduction was applied to $\phi_1\and\psi_1,\dots,\phi_n\and\psi_n$.
\item [$k\equiv$]$2\;\;$(mod 11)

 Let $\phi_1\vee\psi_1,\dots,\phi_n\vee\psi_n$ be all the disjunctions in $\Gamma$ such that reductions have not been applied to them on the previous steps. This case of ($\vee$L) is similar to ($\and$R) reduction.
\item [$k\equiv$]$3\;\;$(mod 11)

 Let $\phi_1\vee\psi_1,\dots,\phi_n\vee\psi_n$ be all the disjunctions in $\Delta$ such that reductions have not been applied to them on the previous steps. This case of ($\vee$R) is similar to ($\and$L) reduction.
\item [$k\equiv$]$4\;\;$(mod 11)

 Let $\phi_1\imp\psi_1,\dots,\phi_n\imp\psi_n$ be all the implications in $\Gamma$ such that reductions have not been applied to them on the previous steps. We add $S(n,2)$ nodes above $\Gamma\Imp\Delta$ and write there all the sequents
$$\psi_{i_1},\dots,\psi_{i_k},\Gamma\Imp\phi_{j_1},\dots,\phi_{j_l},\Delta,$$
where $\{\{i_1\dots i_k\}, \{j_1 \dots j_l\}\}$ is a partition of the set $\{1\dots n)\}$.

We say that a ($\imp$L) reduction was applied to $\phi_1\imp\psi_1,\dots,\phi_n\imp\psi_n$.
\item [$k\equiv$]$5\;\;$(mod 11)

 Let $\phi_1\imp\psi_1,\dots,\phi_n\imp\psi_n$ be all the implications in $\Delta$ such that reductions have not been applied to them on the previous steps. We add a node above $\Gamma\Imp\Delta$ and write there
$$\Gamma\Imp\psi_1,\dots,\psi_n,\Delta.$$
We say that a ($\imp$R) reduction was applied to $\phi_1\imp\psi_1,\dots,\phi_n\imp\psi_n$.

\item [$k\equiv$]$6\;\;$(mod 11)

 Let $\forall x_i\phi_i (x_i)$ be all the formulas of the form $\forall x \phi(x)$ in $\Gamma$.

For every $i$, let $a_i$ be the first available variable, that was not used in a reduction of $\forall x_i\phi_i (x_i)$ before. If there is no such available variable, take first free variable that is not available and mark it as available.

We add a node above $\Gamma\Imp\Delta$ and write there
$$ \phi_1(a_1),\dots,\phi_n(a_n),\Gamma\Imp\Delta. $$
\item [$k\equiv$]$7\;\;$(mod 11)

 Let $\forall x_i\phi_i (x_i)$ be all the formulas of the form $\forall x \phi(x)$ in $\Delta$ such that reductions have not been applied to them on the previous steps. 

Let $a_1\dots a_n$ be the first not available variables. We add a node above $\Gamma\Imp\Delta$ and write there
$$ \Gamma\Imp\phi_1(a_1),\dots,\phi_n(a_n),\Delta. $$
\item [$k\equiv$]$8\;\;$(mod 11)

  Let $\exists x_i\phi_i (x_i)$ be all the formulas of the form $\exists x \phi(x)$ in $\Gamma$ such that reductions have not been applied to them on the previous steps. 

Let $a_1\dots a_n$ be the first not available variables. We add a node above $\Gamma\Imp\Delta$ and write there
$$\phi_1(a_1),\dots,\phi_n(a_n), \Gamma \Imp\Delta. $$
\item [$k\equiv$]$9\;\;$(mod 11)

 Let $\exists x_i\phi_i (x_i)$ be all the formulas of the form $\exists x \phi(x)$ in $\Delta$.

For every $i$, let $a_i$ be the first available variable, that was not used in a reduction of $\exists x_i\phi_i (x_i)$ before. If there is no such available variable, take first free variable that is not available and mark it as available.

We add a node above $\Gamma\Imp\Delta$ and write there
$$ \Gamma\Imp\phi_1(a_1),\dots,\phi_n(a_n), \Delta. $$
\item [$k\equiv$]$10\;\;$(mod 11)

 If $\Gamma$ and $\Delta$ have no common formulas, we add a node above and write there $\Gamma\Imp\Delta$.

\een
It can easily be seen that if every branch of the reduction tree is finite, the sequent is derivable in QGPM.

So if a sequent is not derivable, there is an infinite branch in its reduction tree. Take this branch and denote $\tilde\Gamma$ as the union of all antecedents of the sequents written on the nodes of this branch; $\tilde\Delta$ - the union of all succedents.

Now we define the model: let $M$ be the set of all free variables.

The valuation $v$ is defined by induction on the complexity of formulas. The complexity measures the number of implications:
\ben
\item $c(\phi) = 0$, if $\phi$ is atomic;
\item $c(\phi\and\psi) = c(\phi\vee\psi) = max\{c(\phi),c(\psi) ,\}$;
\item $c(\phi\imp\psi) = c(\psi)+1$;
\item $c(\forall x\phi[x/a]) = c(\exists x\phi[x/a]) = c(\phi(a))$.
\een
Now we define partial n-valuations $v_n$ assigning $0$ and $1$ to atomic formulas and implications of complexity at most n. For every $v_n$, the forcing relation $\model_n\phi$ is defined as usual for all formulas of complexity at most n.

For every n,
\ben
\item $v_n(\phi) = 1 \Iff \phi\in\tilde\Gamma$, if $\phi$ is atomic;
\item $v_n(\phi\imp\psi) = 1\Iff ((\phi\imp\psi)\in\tilde\Gamma\text{ or } \model_{n-1}\psi)$, if $c(\phi\imp\psi)\leqslant n$.
\een
Since partial valuations $v_n$ extend each other, as a limit of $v_n$ we obtain the valuation $v$ (and the relation $\model\phi$) on all closed formulas such that
\ben
\item $v(\phi) = 1 \Iff \phi\in\tilde\Gamma$, if $\phi$ is atomic;
\item $v(\phi\imp\psi) = 1\Iff ((\phi\imp\psi)\in\tilde\Gamma\text{ or } \model\psi)$.
\een
This automatically means that if $\phi\in\Delta$, then $v(\phi)=0$ for atomic $\phi$ and implications (because $\Gamma\cap\Delta = \varnothing$).
\begin{lemma}
For every formula $\phi$,
\ben
\item If $\phi\in\bigcup\Gamma$, then $\model_v \phi$;
\item If $\phi\in\bigcup\Delta$, then $\not\model_v \phi$.
\een
\end{lemma}
\textbf{Proof.}
We prove this by induction on $\phi$.

If $\phi$ is atomic or an implication, we have (1) and (2) by definition.
If $\phi$ is a conjunction or a disjunction, the proof is trivial.

If $\forall x\psi(x)\in\tilde\Gamma$ then for all $a\in M$ $\psi(a)$ was added to $\tilde\Gamma$ on some step. So $\forall a\in M \model \psi(a)$ by induction hypothesis.

If $\forall x\psi(x)\in\tilde\Delta$ than for some $a\in M$ $\psi(a)$ was added to $\tilde\Delta$ on some step. So $\forall a\in M \not\model \psi(a)$ by induction hypothesis.

For $\phi = \exists x\psi(x)$ the proof is similar.
\begin{lemma}
Valuation \emph{$v$} is quasi-boolean.
\end{lemma}
\textbf{Proof.}
We prove that $\model\psi\Imp\model\phi\Imp\psi$ and $\model_v\phi\imp\psi\Imp\not\model_v\phi\text{ or }\model_v\phi$.

Let $\model\psi$. Then $\psi\not\in\tilde\Delta$, so $\phi\imp\psi\not\in\tilde\Delta$. Two cases are possible:
\ben
\item $\phi\imp\psi\in\tilde\Gamma$. Then $v(\phi\imp\psi)=1$ by definition.
\item $\phi\imp\psi\not\in\tilde\Gamma$. Then by definition $v(\phi\imp\psi) = 1$ iff $\model\psi$ iff $\psi\in\tilde\Gamma$. So $v(\phi\imp\psi) = 1$.
\een

Let $\model_v\phi\imp\psi$. Again two cases are possible:
\ben
\item $\phi\imp\psi\in\bigcup\Gamma$. Therefore, $\phi\in\Delta$ or $\psi\in\Gamma$.  If $\psi\in\Gamma$, by the definition $\model_v \psi$. If $\phi\in\Delta$, $\not\model_v\phi$.
\item $\phi\imp\psi\not\in\bigcup\Gamma$. Since $\phi\imp\psi\not\in\bigcup\Delta$, $v(\phi\imp\psi) = v(\psi)$. So $\model_v\psi$.
\een
\qed

Now we can prove completeness. If $\Gamma_0\Imp\Delta_0$ is not derivable in QGPM, take the model built from its reduction tree. Every formula from $\Gamma_0$ is valid in it, every formula from $\Delta_0$ is not valid in it. So this is a counter-model for $\Gamma_0\Imp\Delta_0$. \qed
$$$$

\end{document}